\documentclass[12pt,reqno]{amsart}
\usepackage{amsmath}
\usepackage{amssymb}
\usepackage{a4wide}
\usepackage{soul}
\usepackage{graphicx}
\usepackage{color} 
\usepackage{graphics}
\usepackage{epsfig}
\usepackage{bbm}

\newtheorem{theorem}{Theorem}[section]

\newtheorem{remark}[theorem]{Remark}
\newtheorem{corollary}[theorem]{Corollary}

\newcounter{hypo}

\newenvironment{hyp}{
 \begin{enumerate}
\setcounter{enumi}{\value{hypo}} \item}{\stepcounter{hypo} \end{enumerate}}

\makeatletter
 \@addtoreset{equation}{section}
 \makeatother
 
\title[Numerical approximation of the best decay rate]
{Numerical approximation of the best decay rate for some dissipative systems}

\author[K. Ammari]{Ka\"{\i}s Ammari}
\address{Ka\"{\i}s Ammari, UR Analysis and Control of PDEs, UR 13ES64, Department of Mathematics, 
Faculty of Sciences of Monastir, University of Monastir, 
5019 Monastir, Tunisia}
\email{kais.ammari@fsm.rnu.tn}
\author[C. Castro]{Carlos Castro}
\address{Carlos Castro, Dep. Matem\`atica e Inform\`atica, ETSI Caminos, Canales y Puertos, Universidad Polit\'ecnica de Madrid, 28040
Madrid, Spain}
\email{carlos.castro@upm.es}

\subjclass[2010]{74K10 (35Q72 35B40 34L20)}
\keywords{Rate of decay,  Dissipative systems, Spectral abscissa, Projection method, Eigenvalue approximation.}

\begin{document}

\begin{abstract}
In this paper we propose a numerical method to approximate the best decay rate for some dissipative systems that are bounded perturbation of unbounded skew-adjoint operators. We also give some numerical examples and applications to illustrate the efficiency of this approach.
\end{abstract}

\maketitle

\tableofcontents


\section{Introduction}\label{intro}

Given an evolution dissipative system on a Hilbert space $\mathcal{H}$ with norm $\| \cdot \|$ we define the best decay rate as
\begin{equation} \label{decay_rate}
\inf \{ \omega \mbox{ : } \exists C(\omega)>0 \mbox{ s.t. } \|u(t)\| \leq C \|u(0)\| e^{\omega t}, \mbox{ for every finite energy solution} \}.
\end{equation}  

The characterization of this best decay rate is a
difficult problem and has not a complete answer in the general case. In the particular case of a 
hyperbolic system perturbed with a damping term this question has been addressed by many authors with different techniques 
(see \cite{AmHeTu00_01,AmHeTu01_01,CaCo01_01,CoZu94_01,CoZu95_01,Fr99_01}, and the references therein). These references correspond to one-dimensional problems for which the best decay rate is usually associated to the spectral abscissa of the generator of the
semigroup, i.e. the supremum of the real part of its corresponding eigenvalues. For higher dimension the spectral abscissa also plays a role. In \cite{Le96_01} G. Lebeau characterized the value of the best decay rate for the damped wave equation in terms of two quantities: the spectral abscissa and the mean value of the damping coefficient along the rays of geometrical optics. 

In principle, the spectral abscissa should be easier to compute than the best decay rate, since the latter requires a characterization of the asymptotic behavior of all possible solutions.  However, when considering distributed systems (as wave or beam models) the numerical approximation of the spectrum is not an easy task. For example, if we consider natural approaches based on computing the spectral abscissa of finite elements or finite difference approximations the method fails. This is due to the numerical dispersion introduced by these methods at the high frequencies  that, in particular, affects severely to the approximation of large eigenvalues.  This phenomenon has been largely studied in different contexts specially related with the numerical approximation of control and stabilization problems for hyperbolic models (see for example \cite{z3}, \cite{gll}, \cite{EZ}, \cite{GLH}). Most of the cures proposed so far are based on Tychonoff regularization (see \cite{gll}), filtering techniques for the high frequencies (see \cite{IZ}), multigrid techniques (see \cite{NZ}) or mixed finite element methods (see \cite{CM} and \cite{CMM}). However, none of these techniques seem to provide a uniform approximation of the whole spectrum, required to approximate the spectral abscissa. 

Note however that, even if the whole spectrum of the discrete approximation is not close to the continuous one, we can derive an strategy to obtain a partial approximation, by removing the high frequencies from the finite element model (or the finite differences one). This was the idea followed by M. Asch and G. Lebeau in \cite{ashlebeau} for the 2-D damped wave equation. But, how many frequencies do we have to remove in order to have a good approximation? This is completely unclear in general.  To have an idea, we consider the simplest example of the 1-D undamped wave equation in the interval $x\in(0,1)$. The associated eigenvalues are $\lambda_k= i k\pi$, $k\in \mathbb{Z}^*$, while the finite elements approximation with uniform mesh provides $\lambda_k= i 2N \sin(\frac{k\pi}{2N})$, where $N$ is the number of elements. Assume that we want a $\varepsilon-$approximation of the $k$ frequency, i.e.
$$
\left| i k\pi-i 2N \sin\left(\frac{k\pi}{2N}\right) \right|\leq \varepsilon. 
$$  
Then, using a Taylor expansion of the $\sin(x)$ near $x=0$, we easily obtain an estimate of $N$, 
$$
N\sim  \frac{k^{3/2}\pi^{3/2}}{2 \sqrt{6}\sqrt{\varepsilon}}.
$$
In particular, $N$ grows as $k^{3/2}$. For example, for the lower $20$ frequencies i.e. $|k|\leq 10$, with $\varepsilon=0.1$ we have to compute the eigenvalues of a matrix with size around $(2N)\times (2N)\sim 440\times 440$.   

Here we propose a more efficient method to approximate the spectrum of dissipative systems which are bounded perturbation of unbounded skew-adjoint operators. Roughly speaking, it consists in a projection method on finite dimensional subspaces generated by the first eigenfunctions of the unperturbed operator. As far as we know, the convergence of this projection method in this context has been considered for the first time by J. Osborn to approximate a finite number of eigenvalues when the unperturbed operator is selfadjoint (see \cite{Osborn}). In the present situation, we deal with perturbations of skew-adjoint operators but the result is easily generalized. We show that, under certain hypotheses, this analysis can be adapted to give a uniform approximation of almost the whole spectra and therefore it can be used to approximate the spectral abscissa of the continuous model from the discrete one. In particular, we give an algorithm to approximate this spectral abscissa. The main drawback of the method is that it requires to compute the eigenfunctions of the unperturbed operator. However, this is well-known for one dimensional problems or higher dimensional ones in special domains (intervals or balls), where this method can be easily applied.     

As an application of this method we give some numerical experiments where we show the behavior of the spectral abscissa for different operators under different damping locations.

The rest of the paper is divided as follows: In section \ref{section2} we describe an abstract setting  
for some dissipative systems where, under certain conditions, the decay rate can be described by the spectral abscissa. In section \ref{section3} we describe the numerical projection method and prove the main result stated in Theorem \ref{th_2}, i.e. the uniform convergence of the spectra up to a small number of frequencies.  We also show how this result can be used to approximate the spectral abscissa. In section \ref{section4} we describe a matrix formulation for the numerical method. Finally, in section \ref{section5} we give several examples.

\section{Abstract setting} \label{section2}

Let $A$ be an unbounded operator on  a Hilbert space $H$ with norm $\Vert\cdot\Vert_H$. We assume that $A$ is self-adjoint, positive 
and with compact inverse. We denote its domain by $\mathcal{D}(A) \subset H$. Associated to $A$ we consider $H_{\frac{1}{2}}=\mathcal{D}(A^{\frac{1}{2}})$ the scaled Hilbert space
with the norm $\Vert z \Vert_{\frac{1}{2}}=\Vert A^{\frac{1}{2}} z\Vert_H$, $\forall z\in H_{\frac{1}{2}}$.

\medskip

We set ${\mathcal H}:=H_{\frac{1}{2}} \times H$, endowed with the inner product:
$$
\Big\langle \left[f,g\right],\left[u,v\right]\Big\rangle_{\mathcal{H}} := 
\langle A^{\frac{1}{2}}f,{A^{\frac{1}{2}}u} \rangle_H+ \langle g, {v}\rangle_H, \quad\text{for all} \,\, [f,g], [u,v]\
\text{in}\ {\mathcal{H}},
$$ 
and consider the first order differential equation,
\begin{equation}
\label{Eqd3} 
\left\{
\begin{array}{ll}
\dot{Y}(t) = {\mathcal A}_{{\mathcal B}} Y(t),\\ 
Y(0)=Y_0\in {\mathcal H},
\end{array}
\right.
\end{equation}
where ${\mathcal A}_{{\mathcal B}}:={\mathcal A}_0-{\mathcal B}: 
{\mathcal D}({\mathcal A}_{{{\mathcal B}}})={\mathcal D}({\mathcal A}_0) \subset {\mathcal H} \rightarrow {\mathcal H},$ with
$$
{\mathcal A}_0= \left(
\begin{array}{cc}
\,\, 0  & I \\
- A \,\,  & 0
\end{array}
\right) : {\mathcal D}({\mathcal A}_0) = {\mathcal D}(A)\times H_{\frac{1}{2}} \subset {\mathcal H}  \rightarrow {\mathcal H},$$
and ${\mathcal B} \in {\mathcal L}(\mathcal{H})$ is a bounded operator satisfying 
$$
\Big\langle{\mathcal B} Y,Y\Big\rangle_{\mathcal H} \geq 0, \, \forall \, Y \in \mathcal{H}.
$$

The operator ${\mathcal A}_0$ is skew-adjoint on ${\mathcal H}$ hence 
it generates a strongly continuous group of unitary operators on ${\mathcal H}$, 
denoted by $\big({\bf S}_0(t)\big)_{t \in \mathbb R}$. Since ${\mathcal A}_{{\mathcal B}}$ is dissipative and onto, 
it generates a contraction semi-group on ${\mathcal H}$, denoted by $\big({\bf S}_{{\mathcal B}}(t)\big)_{t \in \mathbb R^+}$. Therefore system (\ref{Eqd3}) is well-posed. Moreover, it is easy to prove that
\begin{equation} \label{ESTEN}
\frac12 \| Y(0) \|_{\mathcal H}^2-\frac12 \| Y(t) \|_{\mathcal H}^2 = \int_0^t \Big\langle {\mathcal B} Y(s),Y(s)\Big\rangle_{\mathcal H} ds\geq 0, \, \forall \, t\geq 0.
\end{equation}
In particular, this means that the mapping $t\mapsto \| Y(t) \|_{\mathcal H}^2$ is not increasing. In many applications it is important to know if this
mapping decays exponentially when $t\to+\infty$, i.e. if the
system (\ref{Eqd3}) is exponentially stable, and if so, what is the minimal rate. This motivates the definition of best decay rate given in (\ref{decay_rate}),

A natural way to characterize this optimal decay  rate is through the spectrum of ${\mathcal A}_{{\mathcal B}}$, denoted by $\sigma({\mathcal A}_{{\mathcal B}})$, since particular solutions associated to a single eigenfunction of ${\mathcal A}_{{\mathcal B}}$ will decay as the real part of the associated eigenvalue. Thus, if we define the spectral abscissa of the generator  ${\mathcal A}_{{\mathcal B}}$ by 
\begin{equation}\label{abscissespec}
\mu(\mathcal{A}_{{\mathcal B}})= \sup \big\{{\rm Re}(\lambda);\ \lambda \in \sigma({\mathcal A}_{{\mathcal B}}) \big\},
\end{equation}
then clearly 
\begin{equation}\label{IN1} 
\mu(\mathcal{A}_{{\mathcal B}})\leq \omega({\mathcal B}).
\end{equation}

Weather this spectral abscissa characterizes the optimal decay rate or not is an interesting question that can be answered in some particular situations. For example, this is true when the eigenfunctions associated to the operator $\mathcal{A}_{{\mathcal B}}$ constitutes a Riesz basis for ${\mathcal H}$, since in this case all solutions can be written as a series of eigenfunctions. 

Let us mention a particular example where this can be proved under certain hypotheses. Consider the damped second order system,
\begin{equation}\label{Eqd1} 
\left\{
\begin{array}{lc}
\ddot{x}(t) + A x(t) + B B^* \dot{x}(t)= 0,\\
\big(x(0),\dot{x}(0)\big)=(x_0,x_1)\in \mathcal{H},
\end{array}
\right.
\end{equation}
with $B$ a bounded operator from $U$ to $H$, 
where $\big(U,\Vert \cdot\Vert_{U}\big)$ is another Hilbert space which will be identified with its dual.

\medskip 

By considering $Y(t)={}^T\big(x(t),\dot{x}(t)\big)$ we can write system (\ref{Eqd1}) in the form (\ref{Eqd3}) with 
$$
{\mathcal B} = \left(
\begin{array}{cc}
0  & \,\, 0 \\
0  & B B^* 
\end{array}
\right)\in {\mathcal L}(\mathcal{H}).
$$

The above system was considered in \cite{AmDiZe13_02}, and a sufficient condition to ensure the existence of Riesz basis constituted by generalized eigenvectors of ${\mathcal A}_{{\mathcal B}}$ was given. The condition concerns the high frequencies of $A$. Let us denote the spectrum of $A$, $\sigma(A)$, by $\{\mu_j\}_{j\geq 1}$ with 
$$
0<\mu_1\leq \mu_2\leq \mu_3 \leq \ldots \leq \mu_n\leq \ldots \to + \infty.
$$ 
For $k\in \mathbb N^*$, we define $\delta_{\pm k}:=\vert \pm i(\sqrt{\mu_{k+1}}-\sqrt{\mu_k})\vert=\sqrt{\mu_{k+1}}-\sqrt{\mu_k}.$ 
We introduce the following assumptions: 

\medskip

\begin{hyp}\label{A1}
$\displaystyle\lim_{k\rightarrow +\infty}\delta_k=+\infty$, 
\end{hyp}

\noindent
and
\begin{hyp}\label{A2}

\,\,\, $\displaystyle \left(\frac{\delta_{k+1}}{\delta_{k}^2}\right)_{k\geq 1}\in l^2(\mathbb N^*)$,
where $ l^2(\mathbb N^*)$ is the space of square integrable sequences.
\end{hyp}

\begin{remark} Observe the following:
\begin{itemize}
\item[(i)] The assumption \ref{A1} implies that the high frequencies of ${\mathcal{A}_0}$ are simple.
\item[(ii)] Assumption \ref{A2} implies
\begin{equation}\label{A3}
\lim_{k\rightarrow +\infty}\left(\frac{\delta_{k+1}}{\delta_{k}^2}\right)=0.
\end{equation}
\item[(iii)] Note that, in general, assumption \ref{A2} does not imply hypothesis \ref{A1}.
\end{itemize}
\end{remark}

\begin{theorem} [Ammari-Dimassi-Zerzeri \cite{AmDiZe13_02}]\label{princb}
Assume \ref{A1} and \ref{A2} hold. Then,
\begin{itemize}
\item[(i)] The eigenvectors of the
associated operator ${\mathcal A}_{\mathcal{B}}$ corresponding to  system (\ref{Eqd1}) form a Riesz basis in the energy space ${\mathcal H}$.
\item[(ii)]
\begin{equation}\label{princss}
\omega({\mathcal B})=\mu(\mathcal{A}_{\mathcal{B}}).
\end{equation} 
\end{itemize}
\end{theorem}
For the abstract-Schr\"odinger equation:
\begin{equation}
\label{scheq}
\left\{
\begin{array}{ll}
\dot{z}(t) + i Az(t) + BB^* z(t) = 0, \forall \, t > 0, \\
z(0) = z_0.
\end{array}
\right.
\end{equation} 
we have the same result as in Theorem \ref{princb} (see \cite{AmDiZe13_02}):
\begin{corollary} \label{schcase}
Assume 

\begin{hyp}\label{A_1}
$\displaystyle\lim_{k\rightarrow +\infty}\delta^\prime_k:= \mu_{k+1}-\mu_k =
+\infty$, 
\end{hyp}

\noindent
and
\begin{hyp}\label{A_2}

\,\,\, $\displaystyle \left(\frac{\delta^\prime_{k+1}}{(\delta^\prime_{k})^2}\right)_{k\geq 1}\in l^2(\mathbb N^*)$.
\end{hyp}

Then,
\begin{itemize}
\item[(i)] The eigenvectors of the
associated operator $A_B := - i A - BB^*$ corresponding to system (\ref{scheq}) form a Riesz basis in the energy space $H$.
\item[(ii)]
\begin{equation}\label{princssh}
\omega(B)=\mu(A_B).
\end{equation} 
\end{itemize}
\end{corollary}

\section{Numerical aproximation of the spectrum} \label{section3}

In this section we consider a projection method to approximate numerically the
spectral abscissa of ${\mathcal A}_{\mathcal{B}}$, i.e. $\mu({\mathcal A}_{\mathcal{B}})$. This projection  method has been used previously to approximate a single eigenvalue, or a fixed number of them, in \cite{Osborn}, when ${\mathcal A}_{\mathcal{B}}$ is a bounded perturbation of a selfadjoint operator.  In the present situation, ${\mathcal A}_{\mathcal{B}}={\mathcal A}_0-{\mathcal B}$ with ${\mathcal A}_0$ skew-adjoint and ${\mathcal B}$ bounded, but the result is easily adapted to this case. We claim that, under certain hypotheses, this analysis provides a uniform approximation of almost the whole spectra and therefore it can be used to approximate the spectral abscissa of the continuous model from the discrete one. In fact, we give below an algorithm to obtain this approximation. 

Along this section we assume that all the eigenvalues  of $A$, that we denote by $\{\mu_k\}_{k\geq 1}$, are simple and ordered increasingly. We also denote by $\{v_k\}_{k\geq 1}$ the associated eigenvectors, that we assume normalized in the norm of $H$. In this case, the eigenvalues of ${\mathcal A}_{0}$ and the corresponding eigenvectors are given by:
\begin{equation}\label{unspectra}
{\mathcal A}_{0}V_{k}=\lambda_k V_{k},\,\,\,{where}\,\, 
  V_{k}=\frac{v_{k}}{\sqrt{2}} \Big[
\frac{1}{\lambda_k}, 1 \Big],\quad \hbox{for all}\quad k\in \mathbb Z^*,
\end{equation}
and $v_{-n}=v_n$, $\lambda_{\pm n}=\pm i\sqrt{\mu_n},$ for $n\in \mathbb N^*$.
Moreover, the family $\big(V_{\pm k}\big)_{k\in \mathbb N^*}$ is an orthonormal basis of the energy
space ${\mathcal H}$.

Let us consider the finite dimensional approximation of ${\mathcal H}$ spanned by the first eigenfunctions of ${\mathcal A}_0$ i.e. 
$$
{\mathcal H}^N = span\{ V_k\}_{k\in Z_N^*}
$$
where $Z_N^*=\{ k\in \mathbb{Z}^*, \; |k|\leq N \}$ and $P^N:{\mathcal H}\to {\mathcal H}^N$ the associated orthogonal projection. 

We define the following Galerkin approximation of the eigenvalue problem (\ref{unspectra}): find $\lambda\in \mathbb{C}$ such that there exists a solution $W^N\in {\mathcal H}^N$, $W^N\neq 0$ of the system 
\begin{equation} \label{unspectraN}
P^N {\mathcal A}_{\mathcal{B}} W^N =\lambda W^N.
\end{equation} 
The spectrum of $P^N {\mathcal A}_{\mathcal{B}}$ is denoted by $\sigma(P^N{\mathcal A}_{\mathcal{B}}) $ and contains $2N$ eigenvalues counting multiplicity.  
Associated to this finite dimensional spectral problem we define the 
spectral abscissa,
\begin{equation}\label{abscissespecN}
\mu(P^N{\mathcal A}_{\mathcal{B}})= \max \big\{{\rm Re}(\lambda);\ \lambda \in \sigma(P^N{\mathcal A}_{\mathcal{B}}) \big\}.
\end{equation}
Our main objective is to relate $\mu(P^N{\mathcal A}_{\mathcal{B}})$ with $\mu({\mathcal A}_{\mathcal{B}})$. Obviously the most we can expect is that $ \mu(P^N{\mathcal A}_{\mathcal{B}})$ approximates the spectral abscissa of the lower $N$ frequencies of ${\mathcal A}_{\mathcal{B}}$. This is basically the statement in Theorem \ref{th_2} below. Note, however, that this is not enough to approximate the spectral abscissa ${\mathcal A}_{\mathcal{B}}$, unless this is given by one of the first $N$ frequencies. 

In practice, there are a number of situations where the large frequencies of ${\mathcal A}_{\mathcal{B}}$ exhibit an asymptotic behavior, in such a way that their real part approaches to a specific value as the frequencies grow. In this case, we only have to consider $N$ sufficiently large to reach this asymptotic regime. Based on this idea we propose below an algorithm to approximate $\mu({\mathcal A}_{\mathcal{B}})$.

Before giving the main result of this section we introduce some notation that is used in the rest of the paper. 
We denote by $\{(\nu_k,U_k)\}_{k\in \mathbb{Z}^*}$ (respectively $\{(\eta_k^N,W_k^N)\}_{k\in \mathbb{Z}^*}$) the eigenvalues and associated eigenvectors of ${\mathcal A}_{\mathcal{B}}$ (respectively $P^N{\mathcal A}_{\mathcal{B}}$), that we assume ordered in such a way that $Im(\nu_{k}) \leq Im(\nu_{k+1})$ and, when equal, $|\nu_{k+1}|\leq |\nu_k|$ (respectively $Im(\eta_{k}^N) \leq Im(\eta_{k+1}^N)$). Note that this is the same ordering chosen for the eigenvalues $\{\lambda_k\}_{k\in \mathbb{Z}^*}$ of $\mathcal{A}_0$. 
We also assume that the above eigenvectors are normalized. 

Finally, to simplify the notation, in this section we omit the space in the norms when there is no confusion, i.e.  we consider $\| \cdot \|=\| \cdot \|_{\mathcal{H}}$ and $\| \cdot \|=\| \cdot \|_{\mathcal{L}(\mathcal{H})}$.

\begin{theorem} \label{th_2}
Assume that the following hypotheses are satisfied:
\begin{itemize}
\item[H1-] ${\mathcal A}_0:D({\mathcal A}_0) \subset \mathcal{H} \to \mathcal{H}$ is skew-adjoint with simple eigenvalues $\{ \lambda_n \}_{n\in \mathbb{Z}^*}$.
\item[H2-] The eigenvalues of ${\mathcal A}_0$ satisfy the following:  
$$
|\lambda_{j+1}-\lambda_j| > 2 \|\mathcal{B}\| , \mbox{ for all $j\in \mathbb{Z}^*$.}
$$ 
\item[H3-] $\|\mathcal{B}\| \leq |\lambda_1|$.
\item[H4-] The eigenvectors of $\mathcal{A}_{\mathcal{B}}=\mathcal{A}_0+\mathcal{B}$ constitute a Riesz basis of $\mathcal{H}$, i.e. there exist constants $m,M>0$ such that
\begin{equation} \label{eq_rb}
m\sum_{j\in \mathbb{Z}^*} |c_j|^2 \leq \left\| \sum_{j\in \mathbb{Z}^*} c_j U_j \right\|^2_{\mathcal H} \leq M\sum_{j\in \mathbb{Z}^*} |c_j|^2, \quad \mbox{ for all $\{ c_j\} \in l^2$}
\end{equation}
\item[H5-] For each $\varepsilon >0$ there exists $r_1>0$ such that
\begin{equation} \label{cond_B}
\max_{|i|\leq p}\sum_{|j|\geq p+r_1} |<\mathcal{B}V_i,V_j>|^2 <\varepsilon, \mbox{for all $p>0$} .
\end{equation} 
\end{itemize}
Then, given $\varepsilon >0$ there exists $r>0$ independent of $N$, such that for all $N>r$
\begin{eqnarray} \label{est_eig1}
\min_j | \eta_p^N-\nu_j|  &\leq& 2\varepsilon, \quad \mbox{ for all } |p|\leq N-r, \\ \label{est_eig2}
 \min_j | \eta_j^N-\nu_p| &\leq& 2\varepsilon, \quad \mbox{ for all } |p|\leq N-r.
\end{eqnarray}
\end{theorem}

\begin{figure}[h]
\centerline{
\includegraphics[height=7cm]{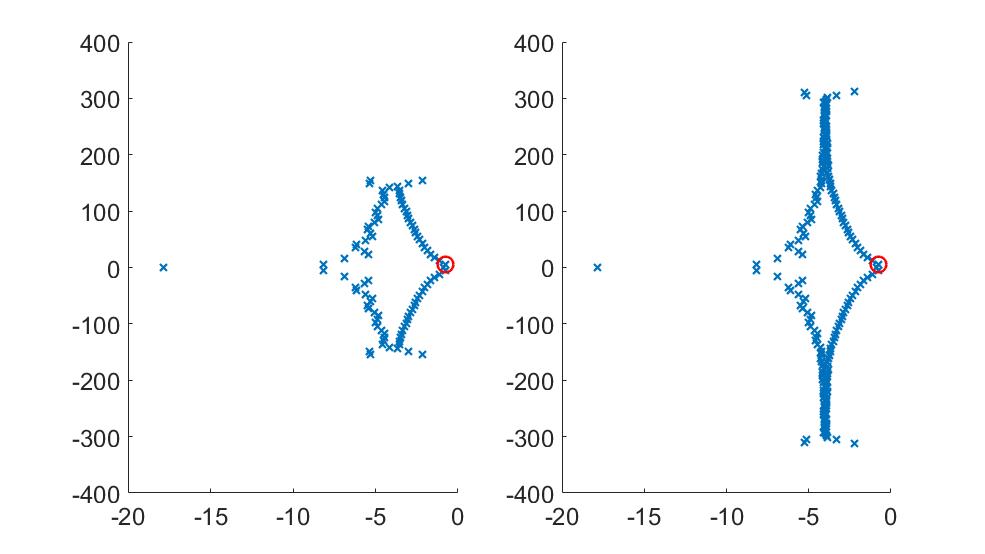} }
\caption{Numerical approximation of the spectrum of the 1-D damped wave equation in the interval $x\in (0,1)$ with a constant damping term $a(x)=10 \chi_{(0.1,0.5)}(x)$, the characteristic function of the interval $(0.1,0.5)$, when considering $N=50$ (left) and $N=100$ (right). The eigenvalue that provides the spectral abscissa is circled. We observe that the highest  frequencies exhibit a clear deviation from the continous ones. However this only happens for a fixed number of frequencies that does not depend on $N$.
\label{fig_00}}
\end{figure}

\begin{remark}
Theorem \ref{th_2} establishes the uniform convergence of the discrete spectrum, as $N\to \infty$, up to the highest $2r$ ones. As we show in the proof below, the value of $r$ in the statement of Theorem \ref{th_2} depends on the value $r_1$ in hypotheses H5 and the assymptotic gap of the eigenvalues of the unperturbed operator, i.e. $|\lambda_{j+1}-\lambda_j|$ for large $|j|$.  Therefore, it can be computed without knowing the eigenvectors of ${\mathcal A}_{\mathcal{B}}$. It can be also estimated numerically for each specific example. In practice (at least in the experiments considered below) this value is small and only a few frequencies must be removed to have uniform convergence. In Figure \ref{fig_00} we show an example for the damping wave equation. In this example $r$ is around $4$ and we observe that this value does not increase for larger values of $N$.  
\end{remark}

\bigskip

\bigskip

\begin{remark}
Hypotheses H2 and H3 in Theorem \ref{th_2} concern the unperturbed operator so that we can check them without knowing the spectrum of ${\mathcal A}_{\mathcal{B}}$. On the other hand, they are used in a technical part of the proof and can be probably relaxed, as checked in the numerical experiments. Hypotheses H4 is more involved since it usually requires some information of the spectrum of ${\mathcal A}_{\mathcal{B}}$. However, the standard technique to prove such property for one-dimensional problems requires basically the asymptotics of the spectrum (see for instance \cite{CoZu94_01}), that is usually known in a number of cases. Finally, hypotheses H5 is a measure of the nondiagonality of the operator ${\mathcal{B}}$ with respect to the basis of eigenvectors for ${\mathcal A}_0$, and can be computed easily. 
\end{remark}

Before giving the proof of Theorem \ref{th_2} we show how this result can be used to approximate $\mu({\mathcal A}_{\mathcal{B}})$ from $\mu(P^N{\mathcal A}_{\mathcal{B}})$.

\begin{corollary} \label{co_1}
Assume that the hypotheses H1-H5 in Theorem \ref{th_2} are satisfied, together with the following one:
\begin{itemize}
\item[H6-] Given $\varepsilon >0$, there exists $N_1$ and $\alpha\in \mathbb{R}$ such that 
\begin{equation}
|\mbox{ Re } (\nu_k) - \alpha| \leq \varepsilon, \quad \mbox{ for all $|k|>N_1$}
\end{equation} 
\end{itemize}
Let us define the following modified spectral abscissa, where we remove the highest frequencies, 
\begin{equation}
\mu_r(P^N{\mathcal A}_{\mathcal{B}})= \max_{|j|\leq N-r} {\rm Re}(\eta^N_j).
\end{equation}
Then, for each $\varepsilon >0$ there exists $N_1$ and $r$, independent of $N$, such that 
\begin{equation} \label{eq_apsp}
|\mu({\mathcal A}_{\mathcal{B}})- \mu_r(P^N{\mathcal A}_{\mathcal{B}})| \leq  \varepsilon, \quad \mbox{ for all $N>N_1$}.
\end{equation}
\end{corollary}

From this result we can easily implement an algorithm to obtain an approximation of the spectral abscissa:
\begin{enumerate}
\item[Step 1.] Choose $\varepsilon>0$ and take $N_1$ from hypotheses H6 above. This value $N_1$ may change from one case to other and we do not know any a priori criteria to compute it. In practice we can give a numerical estimate by chosing $N$ sufficiently large in order to have that the real part of the last frequencies are sufficiently close.
\item[Step 2.] Compute the value $r$ given by Theorem \ref{th_2}. As we have said, this can be computed without knowing the spectrum of ${\mathcal A}_{\mathcal{B}}$. Alternatively, it can be estimated numerically for each specific example too. For example, we choose $N_0$, not too large, and compare the spectra $\sigma(P^{N_0}{\mathcal A}_{\mathcal{B}})$ with $\sigma(P^{2N_0}{\mathcal A}_{\mathcal{B}})$. Then we can estimate $r$ by the smallest value for which 
$$
|\eta_i^{N_0} -\eta_i^{2N_0}| \leq \varepsilon, \quad \mbox{ for all $|i|\leq N-r$} 
$$ 

\item[Step 3.] Take $N>N_1+r$ and compute $\mu_r(P^N{\mathcal A}_{\mathcal{B}})$ by removing the highest frequencies in $\mu(P^N{\mathcal A}_{\mathcal{B}})$. According to Corollary \ref{co_1}, this satifies (\ref{eq_apsp}).  
\end{enumerate}




In the rest of this section we prove Theorem \ref{th_2}.

\bigskip

\begin{proof} (of Theorem \ref{th_2}). We divide the proof in several steps.

{\bf Step 1.} Basic estimates. Here we prove a rough estimate of the eigenvalues. Let 
$$
C_m = \{ \lambda \in \mathbb{C} \mbox{ such that } |\lambda - \lambda_m| <\|\mathcal B \|\}. 
$$
By hypothesis H2, $C_m$ are disjoint. On the other hand, following the perturbation argument in \cite{Osborn}, one can prove that the number of eigenvalues, counting with algebraic multiplicity, of $\mathcal{A}_0, \mathcal{A_{\mathcal{B}}}$ and $ P^N \mathcal{A_{\mathcal{B}}}$ in $C_m$ coincide. In particular, $C_m$ contains only one eigenvalue of $\mathcal{A}_0, \mathcal{A_{\mathcal{B}}}$ and $ P^N \mathcal{A_{\mathcal{B}}}$, i.e.
\begin{equation} \label{est_eig}
|\lambda_k-\nu_k|,\; |\lambda_k-\eta^N_k| < \|\mathcal B \|, \quad \mbox{ for all }k.
\end{equation}

As a consequence, estimates (\ref{est_eig1}) and (\ref{est_eig2}) can be deduced one from the other. We focus on the proof of (\ref{est_eig1}).   

\bigskip

{\bf Step 2.} We prove the following estimate,
\begin{equation} \label{eq_le2a}
\min_{j\in \mathbb{Z}^*} |\nu_j-\eta_p^N|\leq \sqrt{\frac{ M}{m}} \|(I-P^N){\mathcal B} W_p^N\|, \quad  \mbox {for all } |p|\leq N
\end{equation}
where the constants $M$ and $m$ are those given in (\ref{eq_rb}). 

As $\{U_j\}_{j\in \mathbb{Z}^*}$ constitutes a Riesz Basis for ${\mathcal H}$ we can write 
$$
W_p^N=\sum_{j\in \mathbb{Z}^*} \alpha_{p,j} U_j, \quad \alpha_{p,j}\in \mathbb{C}.
$$
We have, 
\begin{eqnarray*}
\sum_{j\in \mathbb{Z}^*} \alpha_{p,j} U_j (\eta_p^N-\nu_j)
&=&P^N{\mathcal A}_{\mathcal{B}} W_p^N-{\mathcal A}_{\mathcal{B}} W_p^N
=P^N({\mathcal A}_0-{\mathcal B})W_p^N -({\mathcal A}_0-{\mathcal B}) W_p^N\\
&=&(I-P^N){\mathcal B} W_p^N.
\end{eqnarray*}
From this identity we have on one hand, 
\begin{equation} \label{eq_le21}
\left\| \sum_{j\in \mathbb{Z}^*} \alpha_{p,j} U^N_j (\eta_p^N-\nu_j)\right\|\leq \left\|(P^N-I){\mathcal B} W_p^N \right\|,
\end{equation}
and on the other hand
\begin{eqnarray*}
&& \!\!\!\! \left\| \sum_{j\in \mathbb{Z}^*} \alpha_{p,j} U^N_j (\eta_p^N-\nu_j)\right\|^2 \geq m \sum_{j\in \mathbb{Z}^*} |\alpha_{p,j}|^2  |\eta_p^N-\nu_j|^2 \geq m \min_{j\in \mathbb{Z}^*}|\eta_p^N-\nu_j|^2 \sum_{j\in \mathbb{Z}^*} |\alpha_{p,j}|^2 \\
&& = \frac{m}{M} \min_{j\in \mathbb{Z}^*}|\eta_p^N-\nu_j|^2,
\end{eqnarray*}
where we have taken into account the normalization of $W^N_p$.
Combining this inequality with (\ref{eq_le21}) we easily obtain 
(\ref{eq_le2a}).

\bigskip

{\bf Step 3.} To estimate the right hand side of (\ref{eq_le2a}) we look at the contribution of the large Fourier coefficients of  $W_p^N$ in the basis $\{ V_j\}_{j\in Z^*}$. In particular we prove that given $r,p_0>0$ such that $ N-r< p_0<N$, the following holds
\begin{equation} \label{eq_3.0}
\| (I-P^{p_0})W_p^N \|\leq \frac{\| \mathcal{B}\|}{|\lambda_{p_0+1}-\lambda_{N-r}|-\| \mathcal{B}\|} , \mbox{ for all $|p|\leq N-r$}. 
\end{equation}

In fact, we have 
\begin{equation} \label{eq_3.1}
\| (I-P^{p_0})W_p^N \|^2 =\sum_{|j|\geq p_0+1} |<W_p^N,V_j>|^2.
\end{equation}
Now observe that
\begin{eqnarray*}
(\eta_p^N-\lambda_j)\Big\langle W^N_p, V_j \Big\rangle&=& \Big\langle P^N {\mathcal A}_{\mathcal{B}} W^N_p, V_j \Big\rangle - 
\Big\langle  W^N_p, ({\mathcal A}_0)^T V_j \Big\rangle = \Big\langle (P^N {\mathcal A}_{\mathcal{B}} -{\mathcal A}_0) W^N_p, V_j \Big\rangle\\
&=&\Big\langle (P^N ({\mathcal A}_0-{\mathcal B})-{\mathcal A}_0) W^N_p, V_j \Big\rangle=-\Big\langle P^N {\mathcal B} W^N_p, V_j \Big\rangle ,
\end{eqnarray*}
and taking modulus in this inequality, 
\begin{equation} \label{eq_3.2}
\left| \Big\langle W^N_p, V_j \Big\rangle \right|^2\leq \frac{\left| \Big\langle P^N {\mathcal B} W^N_p, V_j \Big\rangle \right|^2}{|\eta_p^N-\lambda_j|^2}\leq \frac{\left| \Big\langle P^N {\mathcal B} W^N_p, V_j \Big\rangle \right|^2}{|\eta_{p_0+1}^N-\lambda_{N-r}|^2}\leq \frac{\left| \Big\langle P^N {\mathcal B} W^N_p, V_j \Big\rangle \right|^2}{(|\lambda_{p_0+1}-\lambda_{N-r}|-\| \mathcal{B}\|)^2}, 
\end{equation}
where we have used that $|p|\leq N-r< p_0$, and estimate (\ref{est_eig}) for $\eta_{p_0}^N$, since $N-r\geq m_0$. Therefore, substituting (\ref{eq_3.2}) into (\ref{eq_3.1}), 
\begin{eqnarray*}
\sum_{|j|\geq p_0+1} |<W_p^N,V_j>|^2 &\leq& \sum_{|j|\geq p_0+1}\frac{\left| \Big\langle P^N {\mathcal B} W^N_p, V_j \Big\rangle \right|^2}{(|\lambda_{p_0+1}-\lambda_{N-r}|-\| \mathcal{B}\|)^2}   \\ &\leq& 
\frac{\| P^N\mathcal{B} W^N_p\|^2}{(|\lambda_{p_0+1}-\lambda_{N-r}|-\| \mathcal{B}\|)^2} \leq \frac{\| \mathcal{B}\|^2}{(|\lambda_{p_0+1}-\lambda_{N-r}|-\| \mathcal{B}\|)^2} ,
\end{eqnarray*}
which combined with (\ref{eq_3.1}) gives the desired estimate (\ref{eq_3.0}). 

\bigskip

{\bf Step 4.} Here we obtain (\ref{est_eig1}). First of all, note that for $p_0$ satisfying $N-r>p_0\geq m_0$, we have 
\begin{eqnarray} \nonumber
\|(I-P^N){\mathcal B} W_p^N\| &\leq&  \|(I-P^N){\mathcal B} P^{p_0}W_p^N\| + 
\|(I-P^N){\mathcal B} (I-P^{p_0})W_p^N\| \\ &\leq& \max_{|i|\leq p_0}\|(I-P^N){\mathcal B} V_i\| + \| \mathcal{B} \| \;  \| (I-P^{p_0})W_p^N\| . \label{eq_st4_1}
\end{eqnarray}
Given $\varepsilon >0$, we choose $p_0$ such that the first term in the right hand side of (\ref{eq_st4_1}) is lower than $\varepsilon/2$. From hypothesis $H5$ this can be done for $p_0$ satisfying
\begin{equation} \label{eq_st4_2}
N-p_0\geq r_1, \mbox{ for some $r_1$ independent of $N$}.
\end{equation}
Now we choose $r$ such that the second term in the right hand side of (\ref{eq_st4_1}) is lower than $\varepsilon/2$ too.  From (\ref{eq_3.0}) and hipothesis $H2$ this can be done as long as  
\begin{equation} \label{eq_st4_3}
p_0-(N-r) \geq r_2, \mbox{ for some $r_2$ independent of $N$}. 
\end{equation}
Therefore, if we take $r \geq r_1+r_2$,  estimate (\ref{est_eig1}) follows directly from (\ref{eq_st4_1}). This concludes the proof of Theorem \ref{th_2}.
\end{proof}

\section{Matrix formulation of the numerical method}  \label{section4}     

In this section we reduce the finite dimensional approximation of the eigenvalue problem  (\ref{unspectraN}) to a matrix eigenvalue problem, that we use later to implement the algorithm described in the previous section to approximate the spectral abscissa. 

In order to write the finite dimensional eigenvalue problem (\ref{unspectraN}) in matrix form we write
$$
W^N=\sum_{k\in Z_N^*} c_k V_k
$$
for some Fourier coefficients $c_k$. A straightforward computation shows that if we define 
$$
a_n=\frac{c_n}{\lambda_n}+\frac{c_{-n}}{\lambda_{-n}}, \qquad a_{N+n}={c_n}+{c_{-n}}, \quad n=1,2,...,N
$$
then the eigenvalue problem (\ref{unspectraN}) can be reduced to the equivalent matrix eigenvalue problem 
\begin{equation}\label{eq matrixN}
M_NU^N=\lambda^N U^N, \quad U^N=(a_1,a_2,...a_{2N})^t,
\end{equation}
where $I_N$ is the N-dimensional identity matrix, 
$$
M^N=\left( \begin{array}{ll}  0 & I_N \\ -\Lambda_N & \Omega_N \end{array} \right), \quad
\Lambda_N=\left( \begin{array}{llll}  \mu_1 & 0 &0 & 0\\ 0 & \mu_2& 0 & 0\\ \cdots & \cdots & \cdots & \cdots \\ 0& 0 & 0 & \mu_N \end{array} \right), \quad \Omega_N =\left( \omega_{ij} \right),
$$
and 
$$
\omega_{ij}=-<BB^*v_i,v_j>_H.
$$
The discrete spectral abscissa is defined as 
\begin{equation}
\mu(P^N {\mathcal A}_{\mathcal B}) = \sup \{ \mbox{Re} \lambda^N;\; \lambda_N \in \sigma (M_N) \}.
\end{equation}

\section{Some applications} \label{section5}
We give here some examples of dissipative systems which satisfy or not assumptions \ref{A1} and \ref{A2} but for which we can deduce that the best decay rate can be identified to spectral abscissa. Also, not all the examples satisfy the hypotheses of Theorem \ref{th_2} but the numerical algorithm presented before provides good result in all them.  

\subsection{Damped wave equation} \label{example0}
We consider the following system:
\begin{equation}\label{eq1w}
\partial^2_t u (x,t) - \partial^2_x u(x,t) +
2a(x)\partial_t u(x,t)= 0,\quad
0 < x < 1, \ t > 0,
\end{equation}
\begin{equation}\label{eq2w}
u(0,t) = u(1,t) = 0, \quad t > 0,
\end{equation}
\begin{equation}\label{eq3w}
u(x,0) = u^0(x), \quad \partial_t u(x,0) = u^1(x), \quad
0 < x < 1,
\end{equation}
where $a \in BV(0,1)$ is non-negative satisfying the following condition: 
\begin{equation}\label{condexpw}
\exists \, c>0 \hbox { s.t., } a(x) \geq c,\,\,  \; \hbox{a.e.,\, in a non-empty open subset}  \; I \, \hbox{of}
\;  (0,1).
\end{equation}

We define the energy of the solution $u$ of 
\eqref{eq1w}-\eqref{eq3w}, at time  $t$,  as
\begin{equation}\label{DefEnergyw}
E\big(u(t)\big)=\frac{1}{2}\int_{0}^1 \left( \big|\partial_t u (x,t) \big|^2 +
\big|\partial_x u(x,t)\big|^2\right)\, dx\,.
\end{equation}

$$
U = L^2(0,1), \, H= L^2(0,1), \, H_{\frac{1}{2}} = H^1_0(0,1), $$
$$
{\mathcal D}(A) = H^2(0,1) \cap H^1_0(0,1), 
$$
$$
{\mathcal H} = H^1_0(0,1) \times L^2(0,1), 
$$
$$
A = - \, \frac{d^2}{dx^2}, \quad B \phi = B^* \phi = \sqrt{2a(x)}\phi, \quad \forall \phi \in L^2(0,1).
$$
So,
$$
{\mathcal A}_0 = \left(
\begin{array}{cc}
0 & I \\
 \frac{d^2}{dx^2} & 0
\end{array}
\right), \; {\mathcal A}_{{\mathcal B}} = \left(
\begin{array}{cc}
0 & I \\
 \frac{d^2}{dx^2} & - 2a(x)
\end{array}
\right).
$$

The operator ${\mathcal A}_0$ is skew-adjoint, with compact inverse and
the spectrum is given by $\sigma({{\mathcal A}_0}) = \left\{\pm i k \pi, k \in \mathbb{N}^* \right\}.$ Therefore, the hypotheses (A1) and (A2) are not fulfilled for this problem. However, according to
\cite{Ha89_01}, if $a$ satisfies \eqref{condexpw} then $\omega({\mathcal B}) < 0$. Moreover, in this case it is well-known that the spectral abscissa coincides with the decay rate (see \cite{CoZu94_01}). 

To approximate the spectral abscissa we follow the algorithm described above. This requires to check the hypotheses H1-H6. H1 is clearly true, while H2-H3 depend on the norm of $\mathcal{B}$. In this case, 
$$
\| \mathcal{B} \| \leq 2 \| a\|_{L^\infty (0,1)},
$$
and this is satisfied as soon as $ \| a\|_{L^\infty (0,1)}<\pi/2$. We consider below some examples where this condition is not satisfied but the algorithm works fine. This constitutes a numerical evidence of the non optimality of hypotheses H2-H3.
Hypotheses H4 and H6 were proved for this damped wave equation in  \cite{CoZu94_01} under condition \eqref{condexpw}. Finally, hypotheses H5 has to be checked for each specific example but, in general, it is easy to establish when $a(x)$ is a characteristic function of a subinterval $\omega \subset (0,1)$ or a finite linear combination of the eigenvectors of $\mathcal{A}_0$. 

\begin{figure}[h]
\centerline{\includegraphics[height=5cm]{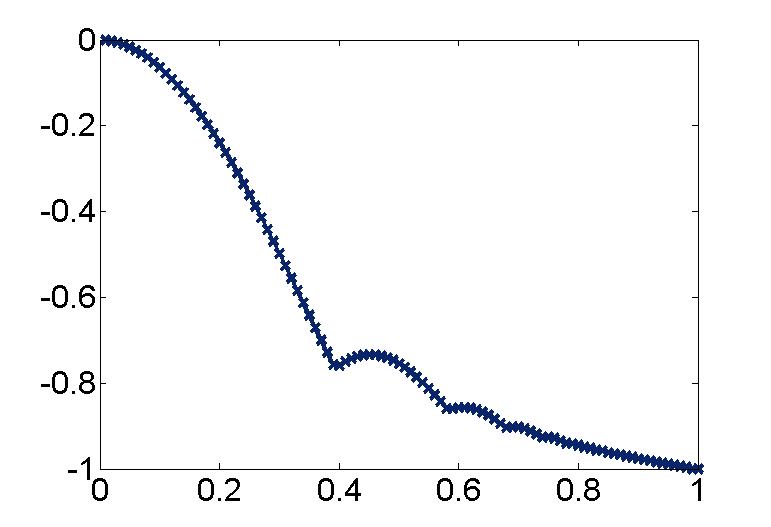} }
\caption{Spectral abscissa versus $\alpha$ for the wave equation when $a(x)$ is given by (\ref{eq_ex1_1}) and $x_0=0.5$.  \label{fig1}}
\end{figure}

The eigenvalue problem is reduced to the matrix eigenvalue problem (\ref{eq matrixN}) for $N=100$ that we solve with MATLAB. In this case, $\mu_k=k^2\pi^2$ and 
\begin{equation} \label{eq_omij}
\omega_{jk}=-2\int_{0}^1 a(x)\sin(j \pi x)\sin(k\pi x) dx.
\end{equation}
In order to see how the distribution of the damping affects to the spectral abscissa  we consider several experiments. We take for $a(x)$ characteristic functions, so that (\ref{eq_omij}) can be computed explicitly, with $\int_0^1a(s)ds=1$ to maintain the same amount of damping. In this way we consider the following two parametric family of dampings,
\begin{equation} \label{eq_ex1_1}
a(x)=\frac1{\alpha} \chi_{(x_0-\alpha/2,x_0+\alpha/2)}(x), \qquad x_0\in (0,1/2), \; \alpha \in (0,2x_0].
\end{equation}
Here $x_0\in (0,1/2]$ is the center of the support and $\alpha\in (0,2x_0]$ its length. For $x_0$ fixed and $\alpha\in (0,2x_0]$ we illustrate the effect of concentrating the damping around a single point. When $\alpha$ approaches zero we formally obtain a Dirac mass concentrated at $x_0$. 
In Figure \ref{fig1} we show the dependence of the spectral abscissa on $\alpha$ when $x_0=1/2$. In particular this spectral abscissa becomes larger as $\alpha \to 0$ and smaller when $\alpha \to 1$ which corresponds to the constant case $a(x)=1$. We also observe that the spectral abscissa is not monotone and that there are some values, around $\alpha=0.4$ and $\alpha=0.6$ for example, where it is not a differentiable function with respect to $\alpha$. In the case $\alpha=0.4$ this corresponds to the case where the spectral abscissa changes from the first to the second eigenvalue. We also observe the presence of several local minima.  

\bigskip

\bigskip

\begin{figure}[h]
\centerline{\includegraphics[height=5cm]{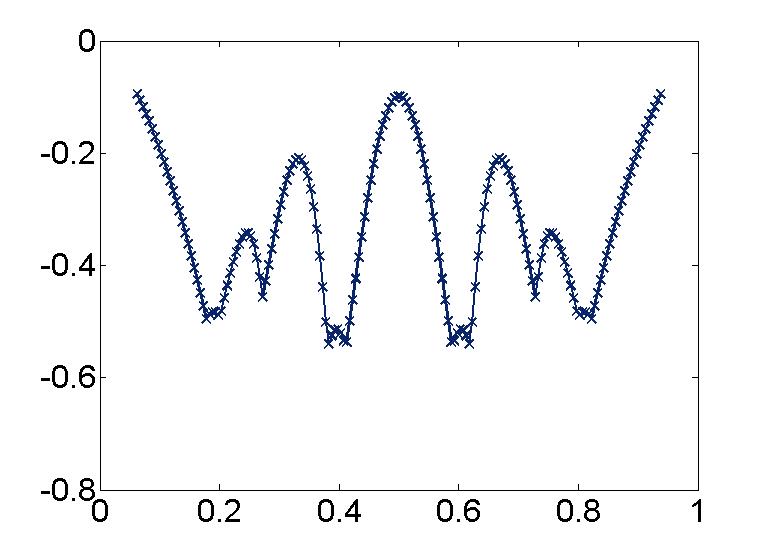} }
\caption{Spectral abscissa versus $x_0$ for the wave equation when when $a(x)$ is given by (\ref{eq_ex1_1}) and $\alpha=1/8$.  \label{fig2}}
\end{figure}

Now, we fix the parameter $\alpha=1/8$ and move the point $x_0\in [1/4,1/2])$. Once again, 
the damping is a one-parametric family of characteristic functions with support in a segment of length $1/8$ that we move through the interval $(0,1)$. The idea is to understand how the location of the damping affects to its efficiency. In Figure \ref{fig2} we show the dependence of the decay on $\alpha$. Note that lower spectral abscissa are obtained around the values $\alpha=0.2,\; 0.4,\; 0.6$ and $0.8$.  

\bigskip

\bigskip

\begin{figure}[h]
\centerline{\includegraphics[height=5cm]{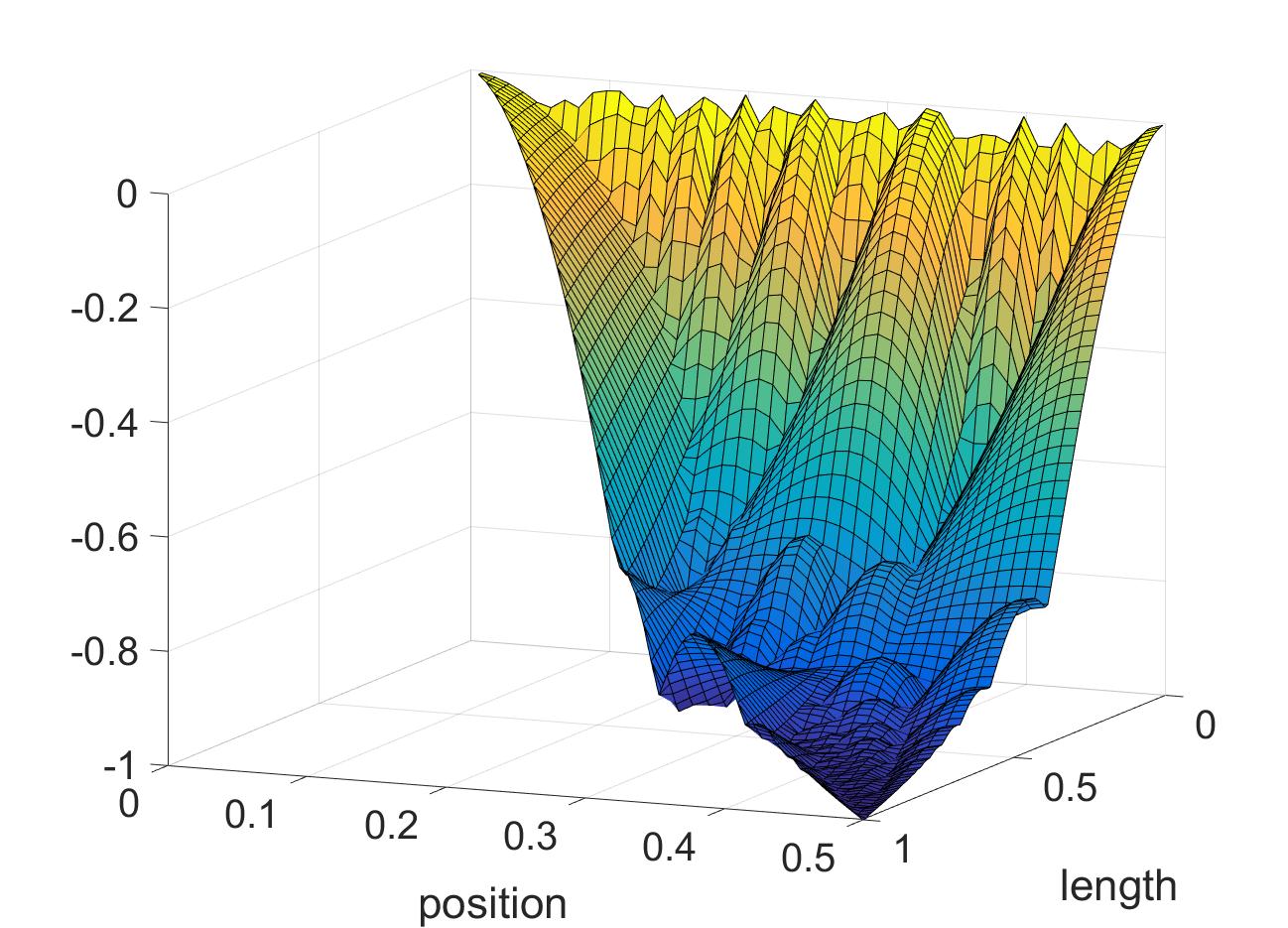} }
\caption{Spectral abscissa versus $(x_0,\alpha)$ for the wave equation when $a(x)$ is given by (\ref{eq_ex1_1}).  \label{fig2_b}}
\end{figure}

In Figure \ref{fig2_b} we show the behavior of the spectral abscissa when we move both parameters. We see that the lower value corresponds to $(x_0,\alpha)=(1/2,1)$ which is when the damping is uniformly distributed in the interval $(0,1)$.

\bigskip

\bigskip

\begin{figure}[h]
\centerline{\includegraphics[height=5cm]{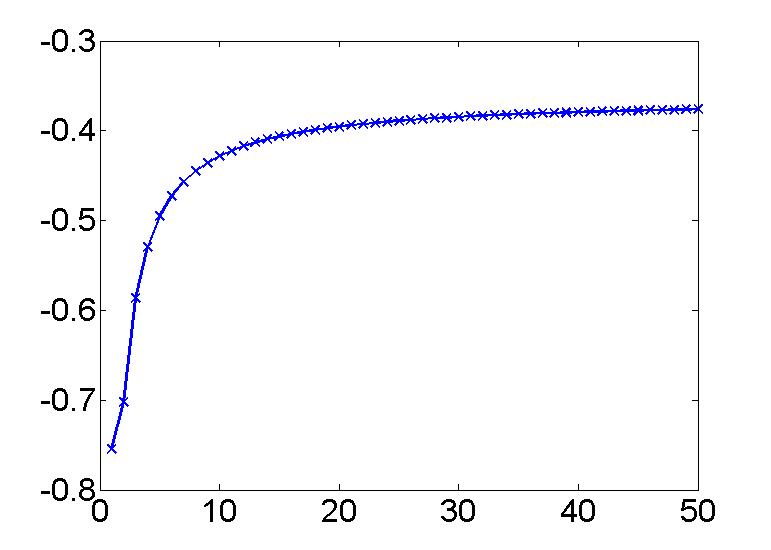} }
\caption{Spectral abscissa versus $\beta$ (number of subientervals where the damping is supported) for the wave equation when $a(x)$ is given by (\ref{eq_ex1_2}).  \label{fig3}}
\end{figure}

Another interesting example appears when the damping is a one-parametric family of a finite number of characteristic functions of the form 
\begin{equation} \label{eq_ex1_2} 
a(x)= \sum_{i=1}^\beta \frac1{2}\chi_{(\frac{2i-1}{2\beta}-\frac1{4i\beta},\frac{2i-1}{2\beta}-\frac1{4i\beta})}(x),\qquad  \beta\in N.
\end{equation}
Note that for larger $\beta$ we consider a damping supported in an increasing number of intervals with lower length in such a way that we still conserve the total mass $1$. In this way we try to understand the influence of an oscillating damping in the spectral abscissa. In Figure \ref{fig3} we show the dependence of the spectral abscissa on $\beta$. We see that as we increase the number of intervals the spectral abscissa becomes larger. Thus, this can be considered as a poor strategy to distribute the damping.  

\bigskip

\subsection{Damped Euler-Bernoulli beam equation} \label{example1}

We consider the following system:
\begin{equation}\label{eq1}
\partial^2_t u (x,t) + \partial^4_x u(x,t) +
2a(x)\partial_t u(x,t)= 0,\quad
0 < x < 1, \ t > 0,
\end{equation}
\begin{equation}\label{eq2}
u(0,t) = u(1,t) = 0, \quad \partial^2_x u(0,t) = \partial^2_x u(1,t) = 0, \quad t > 0,
\end{equation}
\begin{equation}\label{eq3}
u(x,0) = u^0(x), \quad \partial_t u(x,0) = u^1(x), \quad
0 < x < 1,
\end{equation}
where $a \in L^\infty(0,1)$ is non-negative satisfying the following condition: 
\begin{equation}\label{condexpp}
\exists \, c>0 \hbox { s.t., } a(x) \geq c,\,\,  \; \hbox{a.e.,\, in  a non-empty open subset}  \; I \, \hbox{of}
\;  (0,1).
\end{equation}

We define the energy of a solution $u$ of 
\eqref{eq1}-\eqref{eq3}, at time  $t$,  as
\begin{equation}\label{DefEnergyp}
E\big(u(t)\big)=\frac{1}{2}\int_{0}^1 \left( \big|\partial_t u (x,t) \big|^2 +
\big|\partial_x^2 u(x,t)\big|^2\right)\, dx\,.
\end{equation}

$$
U = L^2(0,1), \, H= L^2(0,1), \, H_{\frac{1}{2}} = H^2(0,1) \cap H^1_0(0,1), $$
$$
{\mathcal D}(A) = \left\{u \in H^4(0,1) \cap H^1_0(0,1); \frac{d^2u}{dx^2} (0) = \frac{d^2u}{dx^2} (1) = 0 \right\}, 
$$
$$
{\mathcal H} = [H^2(0,1) \cap H^1_0(0,1)] \times L^2(0,1), 
$$
$$
A = \frac{d^4}{dx^4}, \quad B \phi = B^* \phi = \sqrt{2a(x)}\phi, \quad \forall \phi \in L^2(0,1).
$$
So,
$$
{\mathcal A}_0 = \left(
\begin{array}{cc}
0 & I \\
- \frac{d^4}{dx^4} & 0
\end{array}
\right), \; {\mathcal A}_{{\mathcal B}} = \left(
\begin{array}{cc}
0 & I \\
- \frac{d^4}{dx^4} & - 2a(x)
\end{array}
\right).
$$

The operator ${\mathcal A}_0$ is skew-adjoint and with compact inverse and
the spectrum is given by $\sigma({{\mathcal A}_0}) = \left\{\pm i k^2 \pi^2, k \in \mathbb{N}^* \right\},$ then Assumptions \ref{A1} and \ref{A2} are satisfied. As a direct implication of Theorem \ref{princb}, we have that the best decay rate is given by the spectral abscissa (this result was proved in \cite{AmDiZe13_01}).

To approximate the spectral abscissa we introduce the numerical algorithm described above, where the eigenvalue problem is reduced to the matrix eigenvalue problem (\ref{eq matrixN}). In this case, $\mu_k=-k^4\pi^4$, $\lambda_k=\pm i k^2\pi^2$ and the hypotheses of Theorem \ref{th_2} are satisfied for sufficiently small $a(x)$. More precisely, it is enough to consider $\|a\|_{L^\infty(0,1)} < \pi^2/2$. 

\begin{figure}[h]
\begin{tabular}{cc}
\includegraphics[height=4cm]{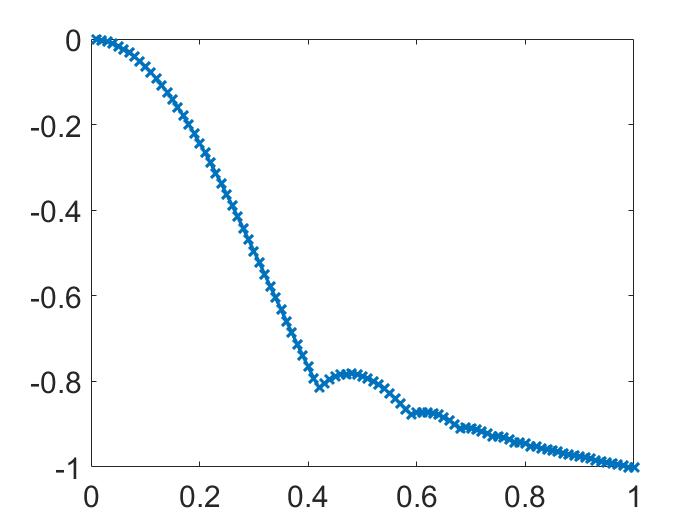} 
& \includegraphics[height=4cm]{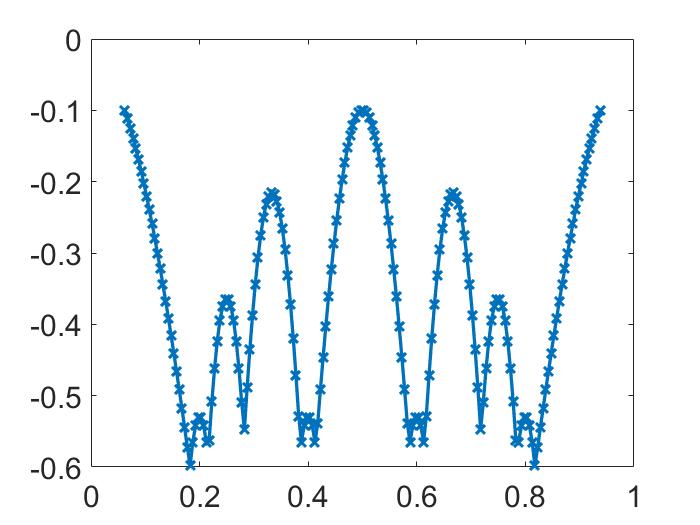} \\
\includegraphics[height=4cm]{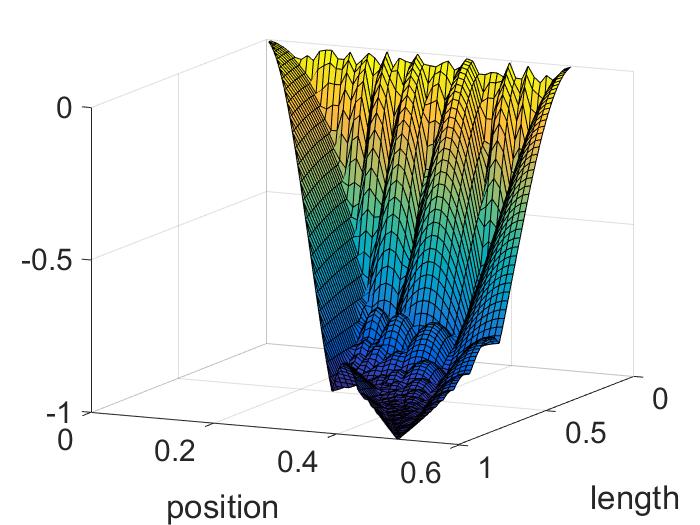} & \includegraphics[height=4cm]{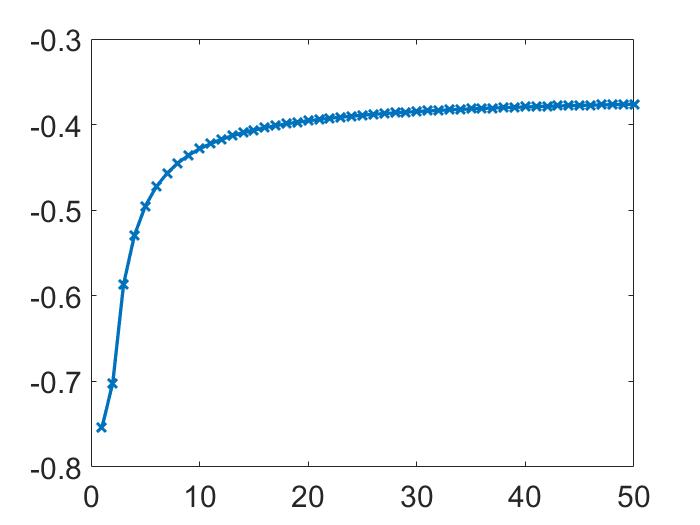}
\end{tabular}

\caption{Spectral abscissa of the beam equation versus $\alpha$ for $x_0=1/2$ (upper left), versus $x_0$ for $\alpha=1/8$ (upper right), versus $(\alpha,x_0)$ (lower left), when $a(x)$ is given by (\ref{eq_ex1_1}), and versus the number of intervals $\beta$ (lower right) when $a(x)$ is given by (\ref{eq_ex1_2}).    \label{fig1_beam}}
\end{figure}

Following the idea of the experiments for the wave equation we have considered the same damping functions $a(x)$. The results are completely analogous as it can be seen in Figure \ref{fig1_beam}.

\bigskip

\subsection{Damped Schr\"odinger equation} \label{example2}

We consider the following system:
\begin{equation}\label{eq1sch}
\partial_t u (x,t) - i \, \partial^2_x u(x,t) +
a(x) u(x,t)= 0,\quad
0 < x < 1, \ t > 0,
\end{equation}
\begin{equation}\label{eq2sch}
u(0,t) = u(1,t) = 0, \quad t > 0,
\end{equation}
\begin{equation}\label{eq3sch}
u(x,0) = u^0(x),  \quad
0 < x < 1,
\end{equation}
where $a \in L^\infty(0,1)$ is non-negative satisfying the following condition: 
$
\exists \, c>0 \hbox { s.t., } a(x) \geq c,\,\,  \; \hbox{a.e.,\, in  a non-empty open subset}  \; I \, \hbox{of}
\;  (0,1).
$

We define the energy of a solution $u$ of 
\eqref{eq1sch}-\eqref{eq3sch}, at time  $t$,  as
\begin{equation}\label{DefEnergy}
E\big(u(t)\big)=\frac{1}{2}\int_{0}^1 \big| u (x,t) \big|^2 \, dx\,.
\end{equation}

In this case, we have 
$$
H =  L^2(0,1), \qquad
A = -  \frac{d^2}{dx^2}, \qquad
B = \sqrt{a(x)} \; I
$$

The operator $i A$ is skew-adjoint and with compact inverse and
the spectrum is given by $\sigma(i A) = \left\{\pm i k^2 \pi^2, k \in \mathbb{N}^* \right\},$ then Assumptions \ref{A_1} and \ref{A_2} are satisfied. As a direct implication of Corollary \ref{schcase}, we have that the best decay rate is given by the spectral abscissa. As in the previous example, the hypotheses of Theorem \ref{th_2} are satisfied for small enough damping terms $a(x)$. Following the idea of the experiments for the wave equation we considered the same damping functions $a(x)$. The results are almost identical to the previous cases (see Figure \ref{fig1_beam}) and we omit them.





\subsection{2D damped wave equation} \label{example3}

We consider the square $\Omega=(0,1)\times (0,1)$ with boundary $\partial \Omega$, and the the following system:
\begin{equation}\label{eq1w2d}
\partial^2_t u (x,t) - \Delta u(x,t) +
2a(x)\partial_t u(x,t)= 0,\quad
x \in \Omega, \ t > 0,
\end{equation}
\begin{equation}\label{eq2wD}
u(x,t) = 0, \quad x\in \partial \Omega, \; t > 0,
\end{equation}
\begin{equation}\label{eq3w2d}
u(x,0) = u^0(x), \quad \partial_t u(x,0) = u^1(x), \quad
x \in \Omega,
\end{equation}
where $a \in BV(\Omega)$ is non-negative satisfying the following condition: 
\begin{equation}\label{condexpw2d}
\exists \, c>0 \hbox { s.t., } a(x) \geq c,\,\,  \; \hbox{a.e.,\, in a non-empty open subset}  \; I \, \hbox{of}
\;  \Omega.
\end{equation}

\medskip

We define the energy of the solution $u$ of 
\eqref{eq1w2d}-\eqref{eq3w2d}, at time  $t$,  as
\begin{equation}\label{DefEnergyw2d}
E\big(u(t)\big)=\frac{1}{2}\int_{0}^1 \left( \big|\partial_t u (x,t) \big|^2 +
\big|\nabla u(x,t)\big|^2\right)\, dx\,.
\end{equation}

$$
U = L^2(\Omega), \, H= L^2(\Omega), \, H_{\frac{1}{2}} = H^1_0(\Omega), $$
$$
{\mathcal D}(A) = H^2(\Omega) \cap H^1_0(\Omega), \quad
{\mathcal H} = H^1_0(\Omega) \times L^2(\Omega), 
$$
$$
A = - \Delta, \quad B \phi = B^* \phi = \sqrt{2a(x)}\phi, \quad \forall \phi \in L^2(\Omega).
$$
So,
$$
{\mathcal A}_0 = \left(
\begin{array}{cc}
0 & I \\
 \Delta & 0
\end{array}
\right), \; {\mathcal A}_{{\mathcal B}} = \left(
\begin{array}{cc}
0 & I \\
 \Delta & - 2a(x)
\end{array}
\right).
$$

\medskip

The operator ${\mathcal A}_0$ is skew-adjoint and with compact inverse and
the spectrum is given by $\sigma({{\mathcal A}_0}) = \left\{\pm i \sqrt{k^2+l^2} \pi, k,l \in \mathbb{N}^* \right\}.$ Note that in this case the properties (A1)-(A2) are not satisfied. However, the spectral abscissa still provides an insight in the decay rate, at least when the support of the damping satisfies the so-called Optics Geometric Condition CGO (see \cite{Le96_01}). 

\begin{figure}[h]
\centerline{\includegraphics[height=7cm]{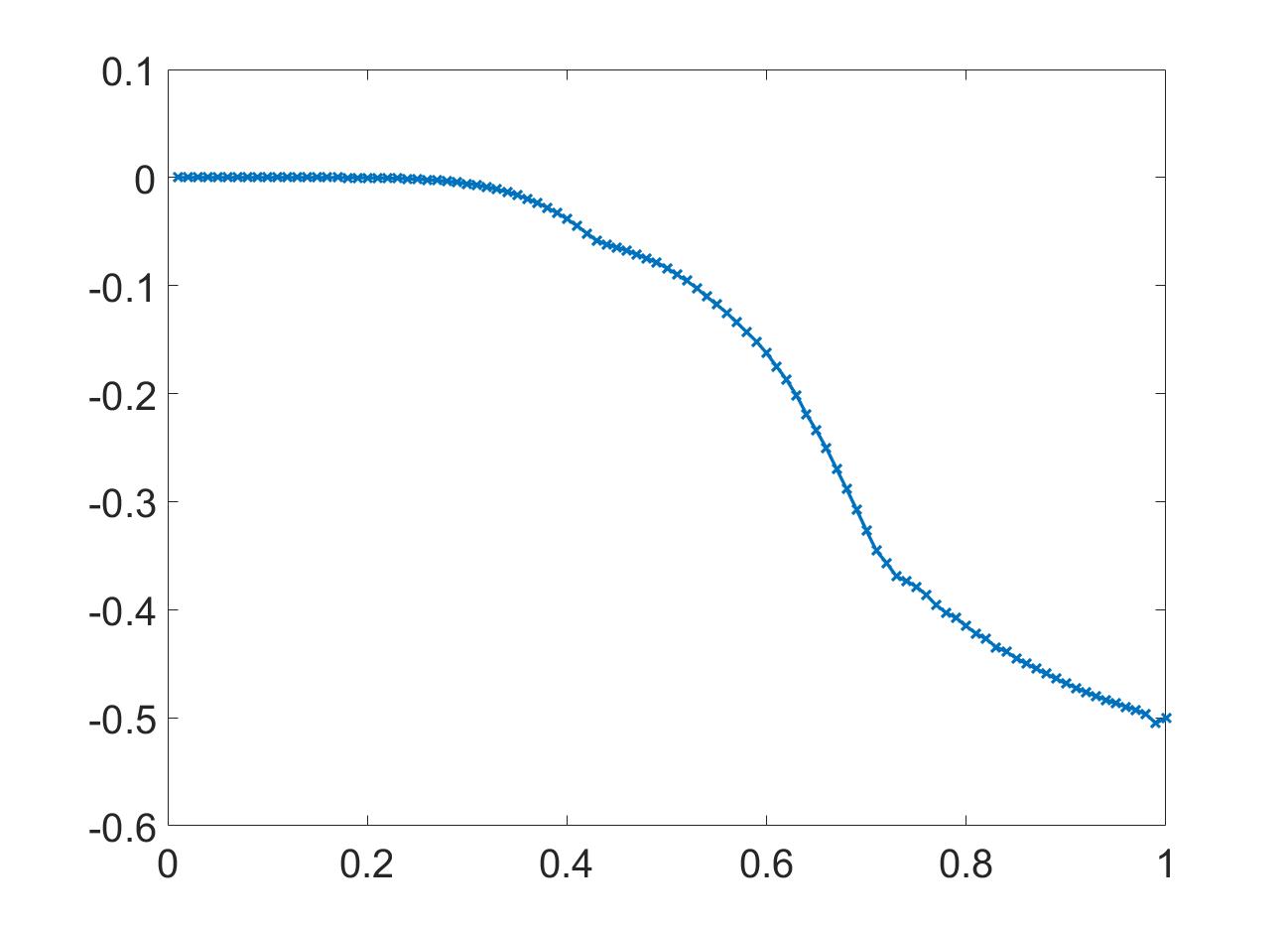} }
\caption{Spectral abscissa versus $\alpha$ for the 2-D damped wave equation when $a(x)=\frac1{\alpha^2} \chi_{D_\alpha}(x)$.  \label{fig1_2dwe}}
\end{figure}

The hypotheses of Theorem \ref{th_2} are not satisfied in this case, but we try the numerical method anyway. Following the idea of the experiments for the 1-D wave equation we analyze the spectral abscissa for damping terms $a(x)$ that are characteristic functions of the form  
\begin{equation} \label{eq_exp_2}
a(x)=\frac1{\alpha^2} \chi_{D_\alpha}(x).
\end{equation}
where $D_\alpha=(0.5-\alpha/2,0.5+\alpha/2)\times (0.5-\alpha/2,0.5+\alpha/2)$, and $\alpha\in(0,1]$. Note that this is a characteristic function of a square that coincides with the domain for $\alpha=1$ and approaches a Dirac as $\alpha\to 0$. 
The results are given in Figure \ref{fig1_2dwe}. Note that in this case the support of the damping does not satisfy the CGO condition and the decay rate at the continuous level can be larger than the spectral abscissa. 

\bigskip

\bigskip

\begin{figure}[h]
\begin{tabular}{cc}
\includegraphics[height=4cm]{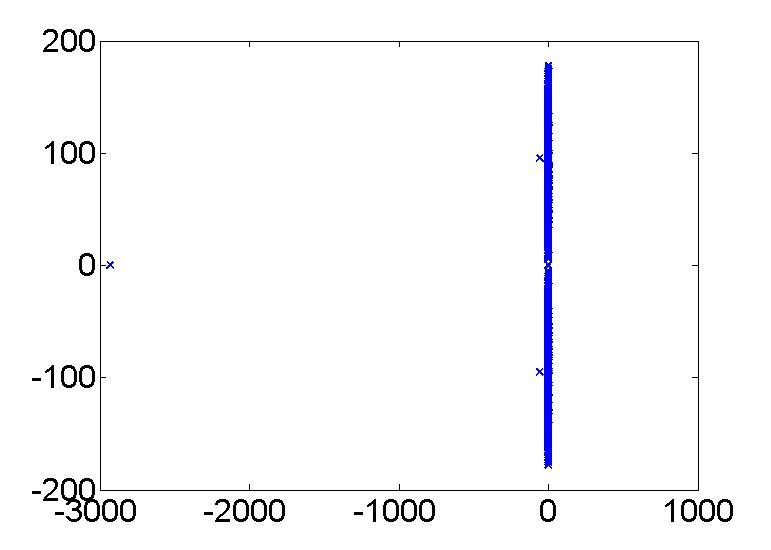} & \includegraphics[height=4cm]{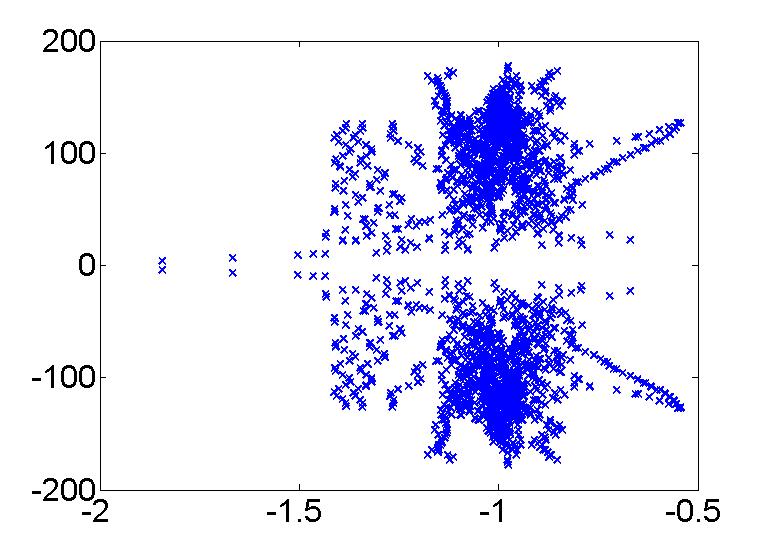} \\
\includegraphics[height=4cm]{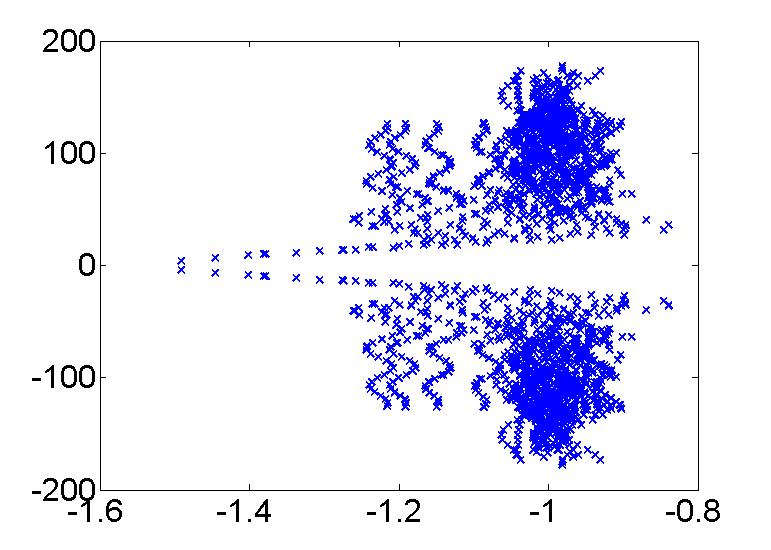} & \includegraphics[height=4cm]{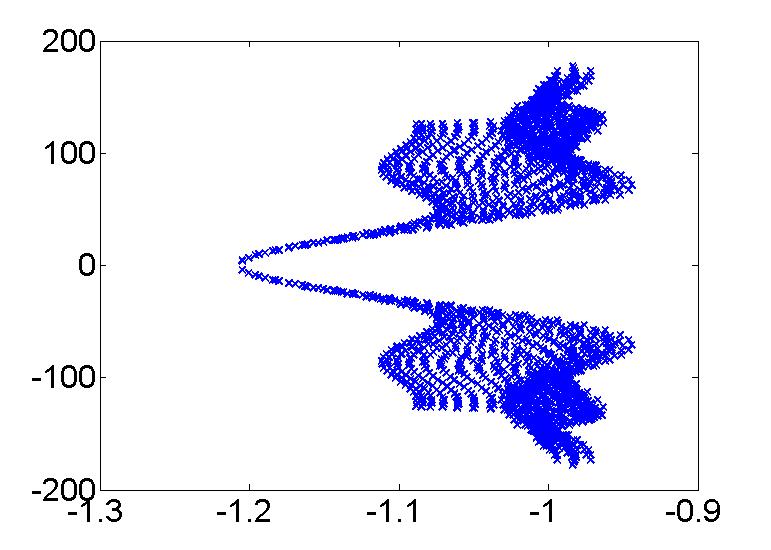} 
\end{tabular}

\caption{Spectrum of the 2-D damped wave equation when $a(x)=\frac1{\alpha^2} \chi_{D_\alpha}(x)$ for different values of $\alpha$ from $\alpha=100$ (upper left) to $\alpha=1$ (lower right).    \label{fig1_2dwe4}}
\end{figure}

We observe that in this case, in contrast with the one-dimensional case, the spectral abscissa is monotone increasing when approaching the Dirac. In Figure \ref{fig1_2dwe4} we show the spectrum in different cases. 

\bigskip

\begin{figure}[h]
\centerline{\includegraphics[height=7cm]{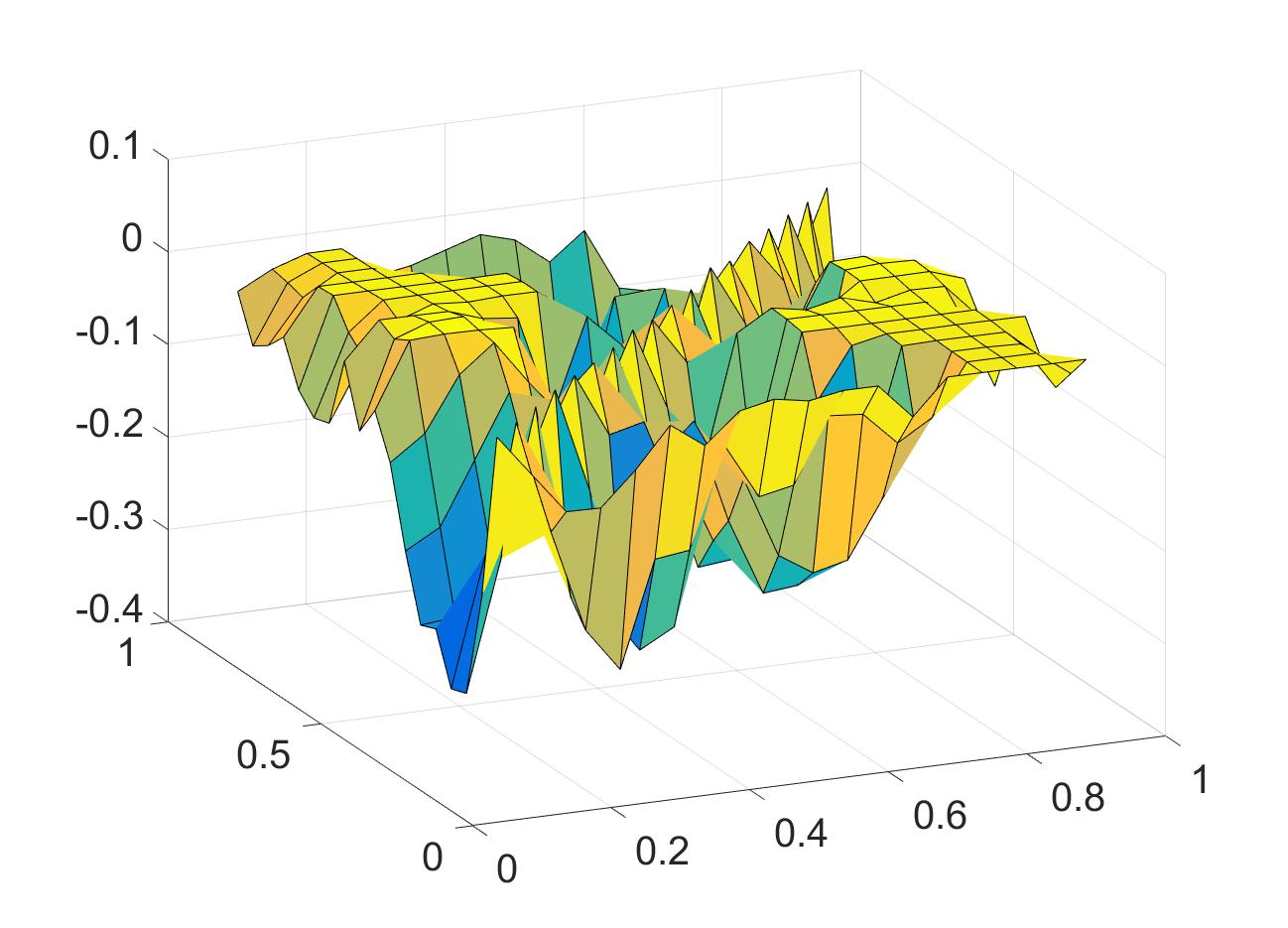} }
\caption{Spectral abscissa versus $\alpha$ for the 2-D damped wave equation when $a(x)= 64\chi_{D_{(\alpha_1,\alpha_2)}}(x) , \qquad \alpha_1,\alpha_2 \in(1/16,15/16)$. Here, $\alpha=(\alpha_1,\alpha_2)$ is the center of the support of the damping.  \label{fig2_2dwe}}
\end{figure}

Now we consider a two-parametric family of the form 
$$
a(x)= 64\chi_{D_{(\alpha_1,\alpha_2)}}(x) , \qquad \alpha_1,\alpha_2 \in(1/16,15/16).
$$ 
where $D_{(\alpha_1,\alpha_2)}$ is the characteristic function of the $\frac18 \times \frac18$ square centered at  $(\alpha_1,\alpha_2)$. Note that now the support of $a(x)$ is a square that we move through the domain $\Omega$, along the two variables, and that we maintain the total mass to $1$. The idea is to understand how the location of the damping affects to the spectral abscissa. In Figure \ref{fig2_2dwe} we show the dependence  on $\alpha$. Note that it shows an oscillating behavior similar to the 1D case.

\section{Conclussions}
We have introduced a projection method to approximate the spectrum of a dissipative system which is a bounded perturbation of a skew-adjoint operator. We show that the associated discrete spectra approximate the frequencies of the continuous problem uniformly with respect to the discretization parameter, up to a fixed number that can be estimated a priori.
Based on this result we introduce an algorithm to approximate the spectral abscissa, and therefore the decay rate, for a large class of dissipative systems. As an applications we analyze the dependence of the damping location in several hyperbolic damped systems.  

\bigskip

{\bf Acknowledgements.}
The second author thanks Prof. E. Zuazua for his comments and suggestions that contributed to improve the presentation of this work.

\end{document}